\documentclass[12pt]{article}
\usepackage{setspace} 
\usepackage[dvips]{graphicx} 
\usepackage{amsmath,amsthm,amsfonts}
\usepackage{epsfig}
\usepackage{subfig}
\RequirePackage[colorlinks,citecolor=blue,urlcolor=blue,linkcolor=blue]{hyperref}
\hypersetup{
colorlinks = true,
citecolor=blue,
urlcolor=blue,
linkcolor=blue,
pdfauthor = {Alexey Kuznetsov},
pdfkeywords = {Hardy-Littlewood function, exponential sum, Euler-Maclaurin formula, Poisson summation},
pdftitle = {Asymptotic approximations to the Hardy-Littlewood function},
pdfpagemode = UseNone
}
    \def\qed{\hfill$\sqcap\kern-8.0pt\hbox{$\sqcup$}$\\}
    \def\beq{\begin{eqnarray}}
    \def\eeq{\end{eqnarray}}
    \def\beqq{\begin{eqnarray*}}
    \def\eeqq{\end{eqnarray*}}

    \def\re{\textnormal {Re}}

    \def\r{{\mathbb R}}

    \def\c{{\mathbb C}}
    \def\Si{{\textnormal {Si}}}

    \def\d{{\textnormal d}}
    \def\i{{\textnormal i}}

	\newtheorem{theorem}{Theorem}
	\newtheorem{lemma}{Lemma}
	\newtheorem{proposition}{Proposition}

	\title{Asymptotic approximations to the Hardy-Littlewood function}
\author{
A. Kuznetsov
\thanks{{Research supported by the
Natural Sciences and Engineering Research Council of Canada.}}  \\ \\
Dept. of Mathematics and Statistics\\  York University \\
4700 Keele Street 
\\Toronto, ON \\ M3J 1P3,  Canada 
 }

\date{Current version: \today}

\begin{document}


\maketitle

\begin{abstract}
\bigskip
The function $Q(x):=\sum_{n\ge 1} (1/n) \sin(x/n)$ was introduced by Hardy and Littlewood \cite{HL1936} in their study of Lambert summability, and since then it has attracted attention of many researchers. 
In particular, this function has made a surprising appearance in the recent disproof by Alzer, Berg and Koumandos \cite{Alzer2005} of a conjecture by Clark and Ismail \cite{Clark2003}. More precisely, Alzer et. al. have shown that the Clark and Ismail conjecture is true if and only if $Q(x)\ge -\pi/2$ for all $x>0$. It is known that  $Q(x)$ is unbounded in the domain $x \in (0,\infty)$ from above and below, which disproves the Clark and Ismail conjecture, and at the same time raises a natural question of whether we can exhibit at least one point $x$ for which $Q(x) < -\pi/2$.  This  turns out to be a surprisingly hard problem, which leads to an interesting and non-trivial question of how to approximate $Q(x)$ for very large values of $x$.  
In this paper we continue the work started by Gautschi in \cite{Gautschi2005} and develop several approximations to $Q(x)$ for large values of $x$. We use these approximations to find an explicit value of $x$ for which $Q(x)<-\pi/2$.  
\end{abstract}

{\vskip 0.5cm}
{\vskip 0.5cm}

\newpage
\section{Introduction}\label{section_intro}

Our main object of interest is the function $Q(x)$, defined as
\beq\label{def_Q}
Q(x):=\sum\limits_{n\ge 1} \frac{1}{n} \sin\left(\frac{x}{n}\right). 
\eeq
It is easy to see that the above series converges for all $x \in \c$ and that $Q(x)$ is an odd entire function of exponential type one. This function was named ``Flett's function" by van de Lune \cite{Lune1981} and D\"oring \cite{Doring1985}, however later it was given the name ``Hardy-Littlewood" function by Alzer et. al. \cite{Alzer2005} and Gautschi \cite{Gautschi2005}, and we will follow the latter convention in our paper.

The function $Q(x)$ was originally introduced in 1936 by 
Hardy and Littlewood \cite{HL1936}, who have used it to construct a certain counter-example in their investigations of Lambert summability. In particular, Hardy and Littlewood have established that $Q(x)=\Omega_+ (\sqrt{\ln\ln(x)})$, which is just a short notation for
\beqq
\liminf_{x\to+\infty} Q(x)/\sqrt{\ln\ln(x)} > 0.
\eeqq
 Hardy and Littlewood also note that the ``O-problems" for $Q(x)$ are much like the corresponding problems for the Riemann zeta function $\zeta(s)$ on the vertical line $\re(s)=1$  (note that it is known that $|\zeta(1+\i t)|=\Omega_{+} ( \ln \ln (t))$, see Theorem 8.5 in \cite{Titchmarsh1986}). 

It seems that this similarity in the behavior of $Q(x)$ and $\zeta(1+\i t)$ is the main reason why $Q(x)$ has attracted so much attention after the work of Hardy and Littlewood. Flett \cite{Flett1950} has proved that $Q(x)=O(\ln(x)^{\frac{3}{4}} (\ln\ln(x))^{\frac{1}{2}+\epsilon})$ and that the same estimate is true for $\zeta(1+\i t)$.
Segal \cite{Segal1972} has discovered  many interesting properties
of $Q(x)$, such as the following identity
\beq\label{identity_Bessel_functions}
\int\limits_0^y Q(x) \d x= 2\sum\limits_{n\ge 1} \sin\left(\frac{y}{2n} \right)^2 = \frac{\pi y}{2}-\frac{1}{2} +  \sum\limits_{m\ge 1} 
\left(\frac{\pi y}{2m}\right)^{\frac{1}{2}} J_1\left( 2 \sqrt{2\pi y m} \right),
\eeq
where $J_1(z)$ is the Bessel function of order one, see \cite{Jeffrey2007}. 
As an application of \eqref{identity_Bessel_functions}, Segal proves that the Ces\`aro, Abel and Borel means of $Q(x)$ are all  equal $\pi/2$. 

Codec\`a \cite{Codeca1984} has studied oscillation and almost periodicity properties of $Q(x)$ and some other related functions.
He also makes an important observation that $Q(x)$ is very similar to some functions which appear as error terms in the estimation of the mean of certain number-theoretic functions. For example, if we define a ``divisor function" $\sigma_{-1}(n):=\sum_{d | n} 1/d$ and $S_{-1}(x):=\sum_{n\le x} \sigma_{-1}(n)$, then it is known (see \cite{Walfisz1936}, p. 100) that
\beq\label{def_S_minus1}
S_{-1}(x)=\frac{\pi^2}{6}x - \frac{1}{2} \ln(x)-\sum\limits_{n\le x} \frac{1}{n} \phi\left(\frac{x}{n} \right)+O(1),
\eeq
where $\phi(x):=\{x\}-\frac{1}{2}$ and $\{x\}$ denotes the fractional part of $x$. Note that both $\sin(x)$ and $\phi(x)$ are periodic functions, which means that the two functions $Q(x)$ and ${\mathcal E}(x):=\sum_{n\le x} (1/n) \phi(x/n)$ should have many similar properties. This turns out to be a beneficial way of looking at these functions: among many other results in \cite{Codeca1984}, Codec\`a has proved that both $Q(x)$  and ${\mathcal E}(x)$ are $B^2$ almost periodic function and are unbounded from above and below. Codec\`a also mentiones that Delange \cite{Delange1980} has extended the results by Hardy and Littlewood 
 \cite{HL1936} and by Flett \cite{Flett1950} and has proved that 
\beq\label{Omega_theorem}
Q(x)=\Omega_{\pm} (\sqrt{\ln\ln(x)}), \;\;\; Q(x)=O((\ln(x))^{\frac{2}{3}}), \;\;\; x\to +\infty.
\eeq
Again, the latter result should be compared with the best known bounds ${\mathcal E}(x)=O((\ln(x))^{\frac{2}{3}})$ (see 
\cite{Walfisz1936}, p. 88) and $\zeta(1+\i t)=O((\ln(t))^{\frac{2}{3}})$ (see formula 6.19.2 in \cite{Titchmarsh1986}). The study initiated by Codec\`a was continued by P\'etermann \cite{Petermann1988}, in particular he gave another proof of \eqref{Omega_theorem}. 

 Van de Lune \cite{Lune1981} has performed the first numerical study of the function $Q(x)$, in particular he has calculated many real and complex zeros of this function. Another method for computing zeros of $Q(x)$ was developed by D\"oring \cite{Doring1985}.

Recently, the Hardy-Littlewood function $Q(x)$ has appeared rather unexpectedly in the dispoof of the Clark and Ismail conjecture by
 Alzer, Berg and Koumandos \cite{Alzer2005}. Clark and Ismail \cite{Clark2003} have studied functions $\Phi_m(x):=-x^m \psi^{(m)}(x)$, where $\psi(x):=\Gamma'(x)/\Gamma(x)$ is the digamma function. They have proved that $\Phi^{(m)}_m(x)$ 
is completely monotone on $(0,\infty)$ for $1\le m \le 16$, and they conjecture that this should be true for all $m \ge 1$. 
However, this conjecture was disproved by Alzer, Berg and Koumandos \cite{Alzer2005} by showing that it is true if and only if $Q(x)>-\frac{\pi}{2}$ for all $x$, the latter statement being false in view of \eqref{Omega_theorem}.
 
Finally, we would like to mention the paper
\cite{Gautschi2005} by Gautschi, which was the main inspiration for our current work. Gautschi has developed two algorithms for computing $Q(x)$ numerically: the first algorithm is based on  the summation by quadrature and the second (more efficient) algorithm is based on truncating the series 
in \eqref{def_Q} at $n = \lfloor x \rfloor$ (the integer part of $x$) and approximating the tail of this series. 
It is clear that this algorithm requires $O(x)$ arithmetic operations in order to obtain a single value of $Q(x)$, thus it becomes impractical if $x$ is very large. 

Our main results in this paper are several new asymptotic approximations for $Q(x)$. Our first result gives an approximation for $Q(x)$, which is extremely accurate in the domain $x>2000$ and which requires only $O(x^{\frac{1}{2}+\epsilon})$ arithmetic operations. The second approximation is somewhat less accurate, but it requires only $O(x^{\frac{1}{3}})$ arithmetic operations. As an application of these two approximations we find that
\beq\label{explicit_x}
Q(8203872394818031742687.4 \times \pi)=-1.5970415\dots,
\eeq
which provides an explicit example to the result by Alzer et. al. \cite{Alzer2005} and answers the question raised by 
Gautschi in \cite{Gautschi2005}.

This paper is organized as follows. In Section \ref{section_results} we review the approximation developed by Gautschi 
\cite{Gautschi2005} and present our main results, Theorems  \ref{thm_big_Oh_of_x_half}, \ref{thm_big_Oh_of_x_third} and
\ref{thm_big_Oh_of_x_eps}. In Section \ref{section_numerics} we perform several numerical experiments: we investigate the accuracy of our approximations, we study the extremes of the function $Q(x)$ and discuss the computations that led to the discovery of \eqref{explicit_x}. 
The proofs of all results are presented in Section \ref{section_proofs}.

\section{Main results}\label{section_results}

Let us introduce the notations that will be used throughout this paper. Given two functions $f: D \mapsto \r$ and $g: D \mapsto (0,\infty)$, we will write $f(x) = O(g(x))$ (or equivalently, $f(x) \ll g(x)$),  if there exists an absolute constant $C>0$ such that $|f(x)|< C g(x)$ for all $x \in D$. When the constant $C$ is not absolute but depends on parameters $p_1,\dots,p_n$, we will write $f(x)=O_{p_1,\dots,p_n}(g(x))$ or $f(x) \ll_{p_1,\dots,p_n} g(x)$. We will also use the notation
$f(x) \approx g(x)$, which stands for $f(x) \ll g(x)$ and $g(x) \ll f(x)$. Bernoulli numbers 
will be denoted by $B_n$. Finally, for $x\in \r$, 
$\lfloor x \rfloor := \max \{ n \le x \; : \; n \in {\mathbb Z}\}$ will denote the integer part of $x$ and $\{x\}=x-\lfloor x \rfloor$ will denote the fractional part of $x$.

Our first result is an algorithm which allows to compute $Q(x)$ to arbitrary precision in $O(x)$ arithmetic operations. 
This algorithm is a simple generalization of the approximation developed by Gautschi (see Section 3 in \cite{Gautschi2005}).

\begin{proposition}\label{prop_big_Oh_of_x}
\emph{($O(x)$ algorithm)}
For any integer $N$ we have
\beq\label{Qx_Ox_algorithm}
Q(x)=\sum\limits_{n=1}^{N-1} \frac{1}{n} \sin\left(\frac{x}{n}\right) + \sum\limits_{k\ge 0} (-1)^k \frac{c_k(N)}{(2k+1)!} \left( \frac{x}{N} \right)^{2k+1},
\eeq 
where $c_k(N)$ are defined as follows 
\beq\label{def_ckN} 
c_k(N):=N^{2k+1} \sum\limits_{n\ge N}  n^{-2k-2}.
\eeq
For large $N$ the coefficients $c_k(N)$ can be computed via the asymptotic expansion
\beq\label{ckN_asymptotic}
c_k(N)=\frac{1}{2k+1}+\frac{1}{2N} + \sum\limits_{m=1}^{M-1} B_{2m} \frac{(2k+2)_{2m-1}}{(2m)!} N^{-2m} + O_{k,M}(N^{-2M}). 
\eeq
\end{proposition}
\vspace{0.25cm}

The proof of the Proposition \ref{prop_big_Oh_of_x} is rather simple, we just sketch the main steps and leave all the details to the reader. In order to obtain \eqref{Qx_Ox_algorithm}, we expand each term of the tail  $\sum_{n\ge N} (1/n) \sin(x/n)$ in the Taylor-Maclaurin series and interchange the order of summation. Asymptotic expansion \eqref{ckN_asymptotic} follows at once by applying the Euler-Maclaurin formula to \eqref{def_ckN}.

As we will see in Section \ref{section_numerics}, for large values of $x$ a good choice of $N$ in Proposition \ref{prop_big_Oh_of_x} is
$N=\lfloor x \rfloor$. Therefore, the above algorithm needs $O(x)$ arithmetic operations and it becomes impractical for very large values of $x$. Our first main result is an asymptotic approximation to $Q(x)$, 
which requires only $O(x^{\frac{1}{2} + \epsilon})$ arithmetic operations and is extremely accurate in the domain $x>2000$. In order to present this algorithm, we will need to define the  sine integral function
\beq\label{def_Si}
\Si(x)=\int\limits_0^x \frac{\sin(t)}{t} \d t. 
\eeq
See Section 8.23 in \cite{Jeffrey2007} for many properties of this function. We will need only one of these properties, namely that
for large values of $x$ the sine integral can be computed via the asymptotic expansion
\beq\label{Si_asymptotic}
\Si(x)=\frac{\pi}{2} - \cos(x)  \sum\limits_{k=0}^{M} (-1)^k \frac{ (2k)!} {x^{2k+1}}
- \sin(x)  \sum\limits_{k=0}^{M}(-1)^k \frac{ (2k+1)!} {x^{2k+2}}+O_M(x^{-2M-3}), \;\; M \in \{-1,0,1,\dots\}.
\eeq
The above asymptotic expansion can be easily derived from formulas 8.215 and 8.233.1 in \cite{Jeffrey2007}.

\begin{theorem}\label{thm_big_Oh_of_x_half}
\emph{($O(x^{\frac{1}{2}+\epsilon})$ algorithm)}
Assume that $0<\epsilon<\frac{1}{2}$. Define $G(x):=\frac{1}{x} \sin\left(\frac{1}{x}\right)$ and $N=N(x):=\lfloor x^{\frac{1}{2}+\epsilon}\rfloor$. Then 
\beq\label{Q_main_sum_1}
Q(x)=\sum\limits_{n=1}^{N-1} \frac{1}{n} \sin\left(\frac{x}{n}\right)
+\frac{1}{2N} \sin\left(\frac{x}{N}\right)+\Si\left(\frac{x}{N}\right)+{\mathcal E}_N(x), 
\eeq
where  ${\mathcal E}_N(x)$ has the following asymptotic expansion: for $M \in \{0,1,2,\dots\}$
\beq\label{error_term}
{\mathcal E}_N(x) =-\sum\limits_{m=1}^{M}  \frac{B_{2m}}{(2m)!} 
\frac{G^{(2m-1)}\left( \frac{N}{x} \right)}{x^{2m}} 
 +O_{M,\epsilon}(x^{-\frac{1}{2}-(4M+3)\epsilon}).
\eeq 
\end{theorem}
\vspace{0.25cm}

The proof of Theorem \ref{thm_big_Oh_of_x_half} will be presented in Section \ref{section_proofs}. This result should be compared with the approximation to the Riemann zeta function $\zeta(s)$ in the critical strip $0 < \re(s) < 1$, see Theorem 4.11 in \cite{Titchmarsh1986}.

Note that in its simplest form the Theorem \ref{thm_big_Oh_of_x_half} gives us the following approximation: for $0<\epsilon<\frac{1}{2}$
\beq\label{Qx_truncated_sum}
Q(x)&=&\sum\limits_{n=1}^{\lfloor x^{\frac{1}{2}+\epsilon}\rfloor} \frac{1}{n} \sin\left(\frac{x}{n}\right)
+\frac{\pi}2 + O_{\epsilon} (x^{-\frac{1}{2}+\epsilon}).
\eeq 
To obtain \eqref{Qx_truncated_sum} one should take $M=-1$ \{  $M=0$ \}  in the asymptotic expansion \eqref{Si_asymptotic}
\{ resp. \eqref{error_term} \} and combine these results with the formula  \eqref{Q_main_sum_1}.
The constant $\pi/2$ in \eqref{Qx_truncated_sum} is reminiscent of the result by Segal \cite{Segal1972} (which was already mentioned in the Introduction), namely that  
the Ces\`aro, Abel and Borel means of $Q(x)$ are all equal $\pi/2$. 
Formula \eqref{Qx_truncated_sum} raises the following two natural questions:
\begin{itemize}
 \item[(i)]  Would it still be correct in the limiting case $\epsilon=0$?
 \item[(ii)] Can we reduce the number of terms in the sum in the right-hand side of \eqref{Qx_truncated_sum}
 to $\lfloor x^{\alpha} \rfloor$ with some $\alpha<\frac{1}{2}$?
 \end{itemize}
 
Our next result provides the answers to both of these questions.

\begin{theorem}\label{thm_big_Oh_of_x_third}
\emph{($O(x^{\frac{1}{3}})$ algorithm)}
Assume that $M$ and $N$ are real numbers, such that $M\ge 1$, $N\ge 1$ and $2\pi M N^2 = x$. Then 
\beq\label{eqn_Qx_third}
Q(x)=\sum\limits_{1\le n <  N} \frac{1}{n}\sin\left(\frac{x}{n}\right)+
\frac{\pi}2+
 \sum\limits_{1\le m < M} 
\left(\frac{\pi}{2mx}\right)^{\frac{1}{4}}
 \sin\left(\frac{\pi}{4}+2\sqrt{2\pi m x} \right)+ O(N^{-1} \ln(x)+x^{-\frac12}N^{\frac{1}{2}}).
\eeq
\end{theorem}
\vspace{0.25cm}

The proof of Theorem \ref{thm_big_Oh_of_x_third} will be presented in Section \ref{section_proofs}. This result should be compared with the approximate functional equation for the Riemann zeta function, see Theorem 4.13 in \cite{Titchmarsh1986}.

Note that if $M>1$ is a fixed constant (which does not depend on $x$), then the second sum in the right-hand side
of \eqref{eqn_Qx_third} is $O(x^{-\frac{1}{4}})$, and the number of terms in the first sum is 
$N=\lfloor  \sqrt{x/(2\pi M)} \rfloor$. 
 This provides the answer to question (i):  we can take $\epsilon=0$ 
 and $N \approx \sqrt{x}$ in \eqref{Qx_truncated_sum}, 
 but then we can only be certain that the error is $O(x^{-\frac{1}{4}} \ln(x))$, and not $O(x^{-\frac{1}{2}})$ as 
\eqref{Qx_truncated_sum} would suggest.
 Another important observation is that the parameters $M$ and $N$ are linked through the condition  $2\pi M N^2 = x$, and 
 if we decrease the number of terms in the first sum in the right-hand side 
 of \eqref{eqn_Qx_third} we will at the same time increase the number of terms in the second sum. 
It is clear that we obtain the best order of approximation and the smallest number of terms in both sums in \eqref{eqn_Qx_third} if we take $M=N=\sqrt[3]{x/(2\pi)}$, in which case the error becomes $O(x^{-\frac{1}{3}} \ln(x))$. In particular, this gives the answer to question (ii). 

\vspace{0.25cm}
\noindent
{\bf Remark 1.}
Segal \cite{Segal1972} has raised the following interesting question. If we formally differentiate the identity \eqref{identity_Bessel_functions}, we would obtain
\beq\label{identity_Bessel_functions2}
Q(x)\stackrel{?}{=}\frac{\pi}{2} + \pi \sum\limits_{m\ge 1} J_0(2\sqrt{2\pi m x}).
\eeq
The problem is that we do not know whether the series in the right-hand side of \eqref{identity_Bessel_functions2} converges.
The connection between formulas 
\eqref{eqn_Qx_third} and  \eqref{identity_Bessel_functions2}  
is provided by the following asymptotic expansion (see formula 8.451.1 in \cite{Jeffrey2007})
\beq\label{J_0_asymptotics}
\pi J_0(2\sqrt{2 \pi m x} )=\left(\frac{\pi}{2mx}\right)^{\frac{1}{4}}
 \sin\left(\frac{\pi}{4}+2\sqrt{2\pi m x} \right) + O\left((mx)^{-\frac{3}{4}} \right),
\eeq
which  implies that \eqref{eqn_Qx_third} could serve as a possible interpretation of \eqref{identity_Bessel_functions2}.
\vspace{0.25cm}

As we have discussed above, the choice $M=N=\sqrt[3]{x/(2\pi)}$ in Theorem \ref{thm_big_Oh_of_x_third} is optimal in the sence that it gives the highest order of approximation and the smallest number of terms in the two sums in the right-hand side of \eqref{eqn_Qx_third}. If we take $N\approx x^{\alpha}$ then the condition 
$2\pi M N^2 = x$ would imply that
$M \approx x ^{1-2\alpha}$, thus is seems that it would be impossible to reduce the number of terms in 
\eqref{Qx_truncated_sum} to $\lfloor x^{\alpha} \rfloor$ with some $\alpha<\frac{1}{3}$. Surprisingly, this is not the case. 
Our next result states that the number of terms  in 
\eqref{Qx_truncated_sum} can be reduced to $\lfloor x^{\alpha} \rfloor$ with  {\it any} $\alpha<\frac{1}{3}$.

\begin{theorem}\label{thm_big_Oh_of_x_eps}
\emph{($O(x^{\epsilon})$ algorithm)}
Let $0<\epsilon<\frac{1}{2}$ and define 
$\delta=\delta(\epsilon):=\epsilon 2^{-\frac{1}{\epsilon}}$.
Then 
\beq\label{eqn_Qx_eps}
Q(x)=\sum\limits_{1\le n < x^{\epsilon}} \frac{1}{n}\sin\left(\frac{x}{n}\right)+
\frac{\pi}2+O_{\epsilon}(x^{-\delta}).
\eeq 
\end{theorem}
\vspace{0.25cm}
The proof of Theorem \ref{thm_big_Oh_of_x_eps} will be presented in Section \ref{section_proofs}. 
To the best of our knowledge, there is no analogous result for the Riemann zeta function (unlike the previous two Theorems). 
\section{Numerical results}\label{section_numerics}

Our first goal is to investigate the accuracy of the approximations provided by 
Theorems \ref{thm_big_Oh_of_x_half} and \ref{thm_big_Oh_of_x_third}. We will use 
Proposition \ref{prop_big_Oh_of_x} in order to compute the benchmark values of $Q(x)$.
We choose $N=\max(\lfloor x \rfloor,1000)$ and compute the coefficients $c_{k}(N)$ for $0\le k \le 9$ using the asymptotic formula \eqref{ckN_asymptotic} with $M=6$. 
Then we truncate the second series in the right-hand side of \eqref{Qx_Ox_algorithm} at $k=k_{\max}=9$. 
Experimenting with higher values of $N$, $M$ and $k_{\max}$ we find that the above parameters guarantee the accuracy of at least 15 digits. 
The code was written in Fortran90 and 
we have used quadruple precision for all computations. 

\begin{figure}
\centering
\subfloat[][$\log_{10}|Q(x)-Q_1(x)|$ for $50\le x \le 5000$]{\label{fig1_error_a}\includegraphics[height =6cm]{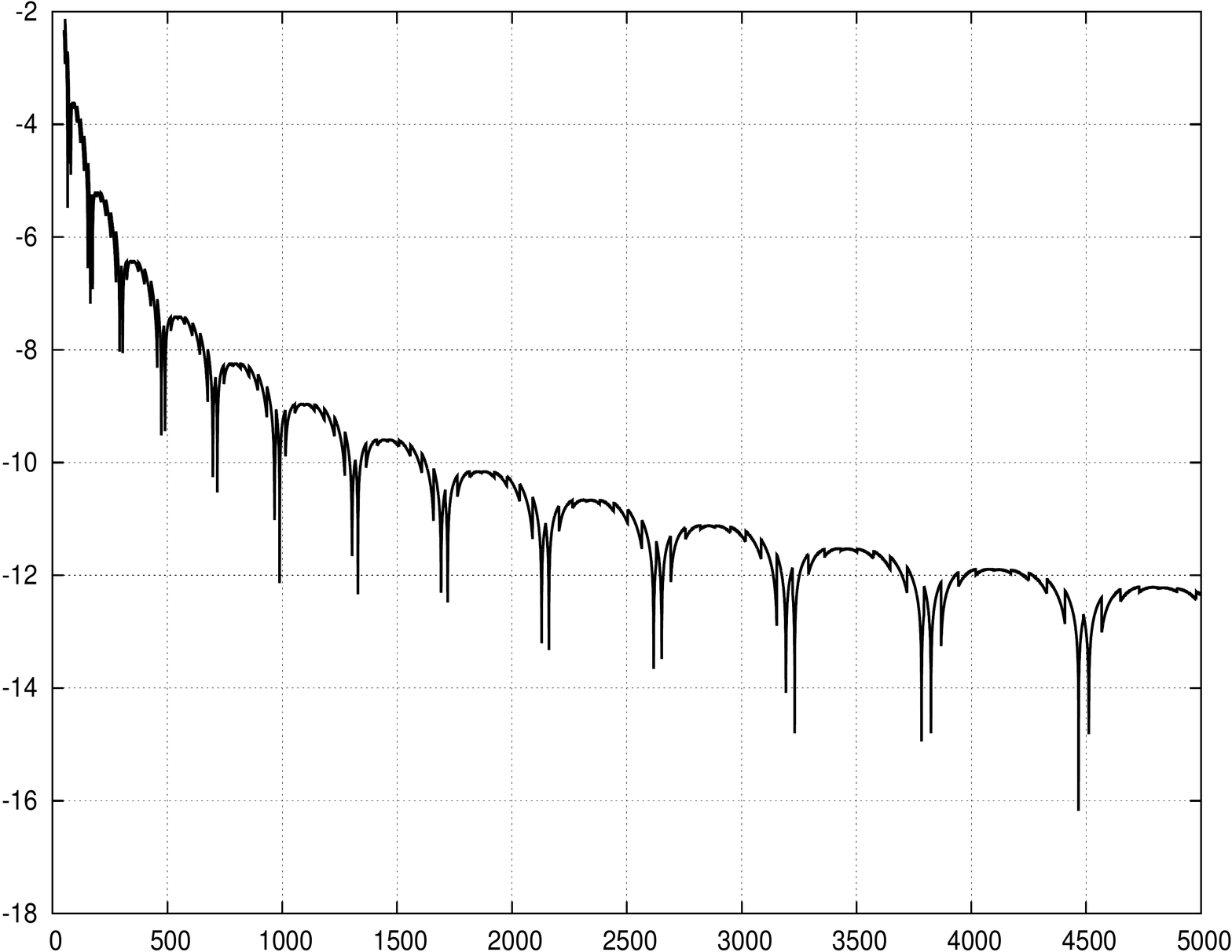}} 
\subfloat[][$Q_2(x_0+x)-Q_1(x_0+x)$ for $0 \le x \le 5\times 10^4$]{\label{fig1_error_b}\includegraphics[height =6cm]{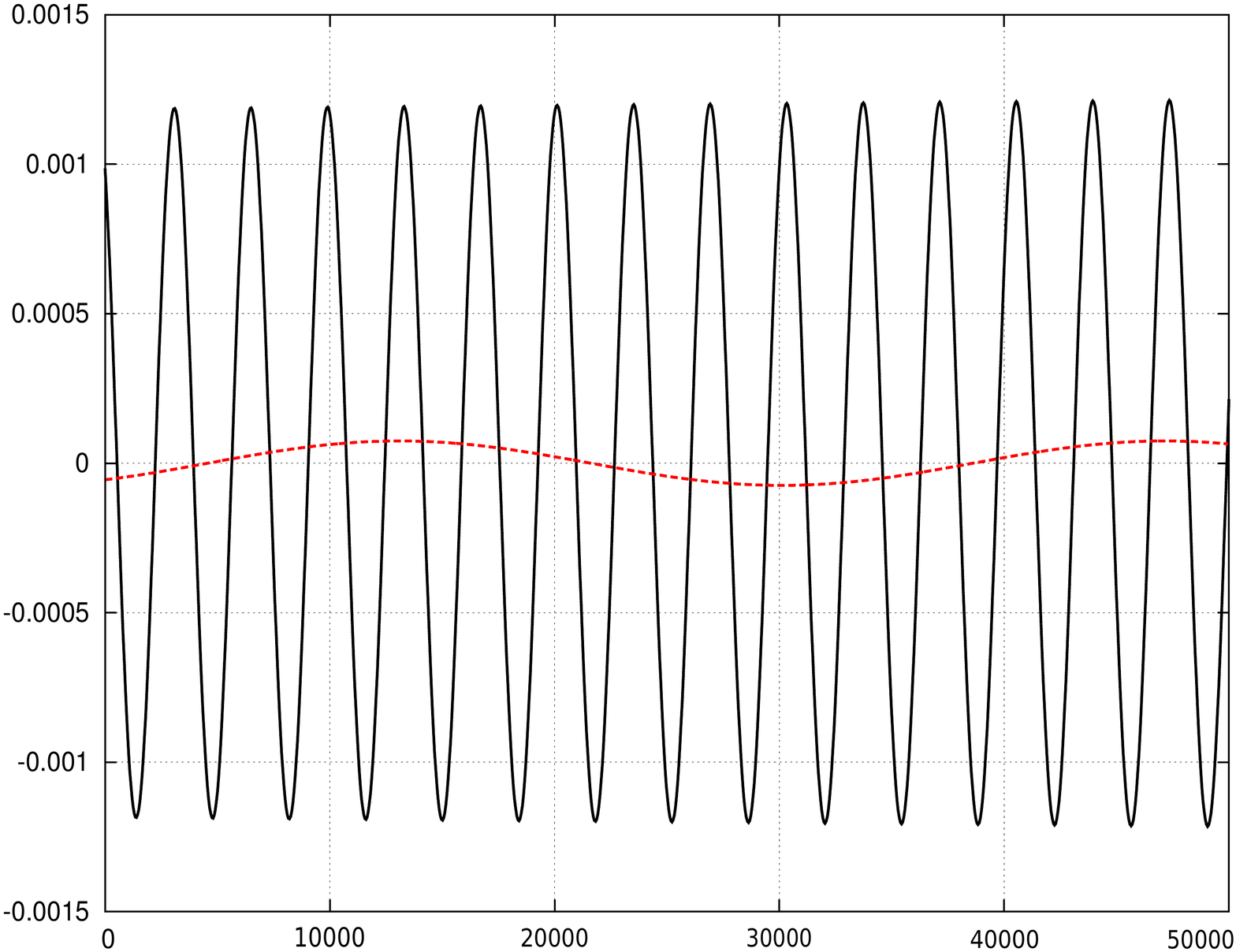}}
\caption{The error of the approximation. Here $Q(x)$ is the value computed using \eqref{Qx_Ox_algorithm}  with $N=\min(x,1000)$ and $M=10$; 
$Q_1(x)$ is computed using \eqref{Q_main_sum_1} and  
\eqref{error_term} with $\epsilon=0.05$ and $M=4$, and $Q_2(x)$ using \eqref{eqn_Qx_third} with $M=N=\sqrt[3]{x/(2\pi)}$. 
  On figure (b) the black curve corresponds to $x_0=10^9$ and the red curve to $x_0=10^{12}$. } 
\label{fig1_error}
\end{figure}

We denote by $Q_1(x)$ the approximation to $Q(x)$ obtained by 
setting $\epsilon=0.05$ and $M=4$ in Theorem \ref{thm_big_Oh_of_x_half}, and by $Q_2(x)$ the approximation to $Q(x)$ 
provided by taking  the optimal choice of parameters $M=N=\sqrt[3]{x/(2\pi)}$ in Theorem \ref{thm_big_Oh_of_x_third}.
We present the results of our computations on figure \ref{fig1_error}.  
On figure \ref{fig1_error_a} we see that $Q_1(x)$ provides an excellent approximation to $Q(x)$: in the domain $x>2000$ the absolute error is  smaller than $10^{-10}$ and 
in the domain $x>4500$ the absolute error is smaller than $10^{-12}$. In order to investigate the accuracy of the second approximation we 
would like to consider much larger values of $x$. Given the fact that for $x$ large  $Q_1(x)$ is very close to $Q(x)$ 
and is much easier to evaluate numerically, we will use the values of $Q_1(x)$ as a benchmark. 
The results are presented on figure \ref{fig1_error_b}. 
We find that the error $Q_2(x)-Q_1(x)$ is of the order $10^{-3}$ when $x\approx 10^9$ and of the order $10^{-4}$ when $x\approx 10^{12}$. Note that these results are in perfect agreement with the error estimate $O(x^{-\frac{1}{3}} \ln(x))$ in \eqref{eqn_Qx_third}. 

As our next goal, we have tried to find a value of $x$ where $Q(x)< -\pi/2$. This would provide an explicit example 
for the key step in the disproof of the Clark and Ismail conjecture by Alzer et. al.  \cite{Alzer2005}.
 Our first attempt was to look at the very large local maximums/minimums of $Q(x)$. 
As was noted by Alzer et. al. \cite{Alzer2005}, for every $n\in \{0,1,2,\dots\}$ the function $Q(x)$ has a local maximum in the interval $[2n\pi, (2n+1)\pi]$ and a local minimum in the interval $[(2n+1)\pi,(2n+2)\pi]$. By 
checking all consecutive local maximums/minimums of $Q(x)$ we have found several values of $x$, where $Q(x)$ has a local extremum and
the value of this local maximum/minimum is greater/smaller than the value of all other local maximums/minimums in the interval $[0,x)$.  
The results of our computations are presented in Table \ref{tab1}. The lowest value of $Q(x)$ that we were able to find using this approach 
is approximately $-1.31$. Given the fact that the function $Q(x)$ grows very slowly (see 
\eqref{Omega_theorem}), it became clear that this brute force search method will not work. Therefore,  we will have to bring in some other ideas in order to reach the value of $Q(x)$ which is less than $-\pi/2$.

In order to find these new ideas, we have looked at the proof of the estimate $Q(x)=\Omega_{\pm} (\sqrt{\ln\ln(x)})$ given in \cite{Alzer2005} and \cite{HL1936}. This proof is based on a rather interesting argument, which is both constructive and non-constructive. 
Let us summarize the main steps of this argument. We define a set of integer numbers ${\mathcal M}$, such that $q\in {\mathcal M}$ if and only if all divisors of $q$ satisfy $q \equiv 1 ({\textnormal{mod}} \; 4)$. The first few elements of ${\mathcal M}$ are 
\beqq
{\mathcal M}=\{5, 13, 17, 25, 29, 37, 41, 53, 61, 65, 73, \dots\}.
\eeqq
Now, for each $k\ge 1$ we define an integer number
\beqq
K=K(k):=\prod\limits_{\substack{1\le q \le 4k+1 \\ q\in {\mathcal M}}} q,
\eeqq
and for each $j\in \{1,2,\dots,K\}$ we define 
\beq\label{def_xkj}
x^{(k)}_j:=(4j+1)K \frac{\pi}{2}, \;\;\; \hat x^{(k)}_j:=(4j+3)K \frac{\pi}{2}.
\eeq
In \cite{Alzer2005} and \cite{HL1936} it was established that there exist constants $a\in \r$, $\hat a \in \r$ and $b>0$, $\hat b >0$, such that for every $k$ there exist $i=i(k)$ and $j=j(k)$ such that
\beq\label{exteme_value_of_Q}
Q(x^{(k)}_i)>a+b \sqrt{\ln(k)}, \;\;\; Q(\hat x^{(k)}_{i})<\hat a-b \sqrt{\ln(k)}.
\eeq
Using this result coupled with the trivial estimate $K(k) < (4k+1)^{4k+1}$ it is not hard to establish that $Q(x) = \Omega_{\pm} ( \sqrt{\ln\ln(x)} )$ as $x\to +\infty$.

\begin{table}
\centering
\bigskip
\begin{tabular}{|c||c|c|c|c|c|}
\hline  \rule[0pt]{0pt}{11pt}
    $x/\pi$  &  12127812.568     & 37324872.600 & 50774112.608 & 176438112.602 & 413620512.598 \\ [0.5ex] \hline 
\rule[0pt]{0pt}{11pt}           $Q(x)$ & 4.3300 &    4.3426 & 4.3931 & 4.4596 & 4.4638   \\ [0.5ex]\hline \hline 
\rule[0pt]{0pt}{11pt}
    $x/\pi$         & 3596987.431  & 3841175.419 & 51836087.411 & 196661495.417 &  580973087.402 \\ [0.5ex] \hline 
\rule[0pt]{0pt}{11pt}           $Q(x)$    & -1.0512   & -1.1635 &  -1.2406  & -1.2740 & -1.3134   \\ [0.5ex]\hline 
\end{tabular}
\vspace{0.2cm}
\caption{Some ``extreme" local maximums/minimums of $Q(x)$}\label{tab1}
\end{table}

\begin{table}
\centering
\bigskip
\begin{tabular}{|c|c|c|c|c|c|}
\hline \rule[0pt]{0pt}{11pt}
    $k$  &   \; $K(k)$ \; &  $\max\limits_{1\le i \le K} Q(x^{(k)}_i)$   &  $\max$ achieved at  & $\min\limits_{1\le j \le K} Q(\hat x^{(k)}_j)$  &  $\min$ achieved at     \\ [0.5ex] \hline \hline
   \rule[0pt]{0pt}{11pt}    4  &  1105 &  4.1352  & $i=876$ & -0.9262 & $j=1103$   \\ [0.5ex]\hline
   \rule[0pt]{0pt}{11pt}    6  &  27625 &  4.2606  & $i=24966$ & -1.1647 & $j=6785$     \\ [0.5ex]\hline
   \rule[0pt]{0pt}{11pt}    7  &  801125 &  4.4127  & $i=259627$ & -1.2498 & $j=468148$     \\ [0.5ex]\hline
   \rule[0pt]{0pt}{11pt}    9  &  29641625 &  4.5752 & $i=7030990$ & -1.4347 & $j=615949$     \\ [0.5ex]\hline
   \rule[0pt]{0pt}{11pt}    10  &  1215306625 &  4.6586 & $i=96330096$ & -1.4717 & $j=69224831$     \\ [0.5ex]\hline
\end{tabular}
\vspace{0.2cm}
\caption{Computing maximum/minimum values of $Q(x^{(k)}_i)$ and $Q(\hat x^{(k)}_j)$ }\label{tab2}
\end{table}

The intuition behind the definition \eqref{def_xkj}  is rather simple: one can see that $\sin(x^{(k)}_i/n)=1$ 
and $\sin(\hat x^{(k)}_j/n)=-1$ for all $n$ which divide $K$. Since $K$ has many small divisors $n$, this shows that there will
be quite many terms in the sum \eqref{def_Q} defining $Q(x^{(k)}_i)$, which are equal to $1/n$ (or which are equal to $-1/n$, in the case of $Q(\hat x^{(k)}_j)$).  Of course, it still takes a lot of work to deduce \eqref{exteme_value_of_Q}, as we have to show that the sum over all $n$ which do not divide $K$ is not too large. See \cite{Alzer2005} and \cite{HL1936} for all the details.

Formulas  \eqref{def_xkj} and \eqref{exteme_value_of_Q} give us an algorithm for finding very large positive/negative values of 
$Q(x)$: these extreme values will happen at points $x^{(k)}_i$ and $\hat x^{(k)}_j$. 
This is the constructive side of this result. 
At the same time, there is no information on how to find the indices $i$ and $j$  for which $Q(x^{(k)}_i)$ or $Q(\hat x^{(k)}_j)$
would achieve the maximum/minimum values, therefore this part would still have to be done by brute force search, by checking 
all the indices in the range $1\le i \le K(k)$. Clearly, this is only feasible if $k$ is not too large, as $K(k)$ grows very fast as $k$ increases.

The results of our computations are presented in Table \ref{tab2}. For $k=4$ and $k=6$ we have used the approximation $Q_1(x)$, provided by Theorem 
\ref{thm_big_Oh_of_x_half} with $\epsilon=0.025$ and $M=4$. 
For larger values of $k$ we have used the approximation $Q_2(x)$, obtained from Theorem 
\ref{thm_big_Oh_of_x_third} by setting $M=N=\sqrt[3]{x/(2\pi)}$. In this case the value of $Q(x)$ at the extreme point was confirmed by computing it via the more accurate approximation $Q_1(x)$. All computations were performed on a regular desktop computer with Intel i7 2600 quad-core processor, running Ubuntu Linux. 
In order to fully utilize all four cores of the processor, we have parallelized the Fortran90 code using OpenMP API.

It is instructive to compare the results presented in Tables \ref{tab1} and \ref{tab2}. 
First of all, one can check that the points where $Q(x)$ attains extremely large local maximums/minimums in Table \ref{tab1} are not located near the points $x^{(k)}_i$ or $\hat x^{(k)}_j$ defined by \eqref{def_xkj}. 
Second, we see that the large values of $Q(x)$ in Table \ref{tab1} happen at smaller values of $x$ than those 
in Table \ref{tab2}. This provides a compelling numerical evidence that the result 
$Q(x) = \Omega_{\pm} ( \sqrt{\ln\ln(x)} )$ (which is derived via \eqref{exteme_value_of_Q})
is suboptimal. A possible way to prove a stronger result would be to  understand how to predict the points in Table \ref{tab1}, where $Q(x)$ attains its extreme local
maximums/minimums. Given the similarity between $Q(x)$, $\zeta(1+\i t)$, ${\mathcal E}(x)$ and other functions arising in Number Theory, this understanding could potentially lead to establishing stronger results about oscillation of all these functions.

\begin{figure}
\centering
\includegraphics[height =8cm]{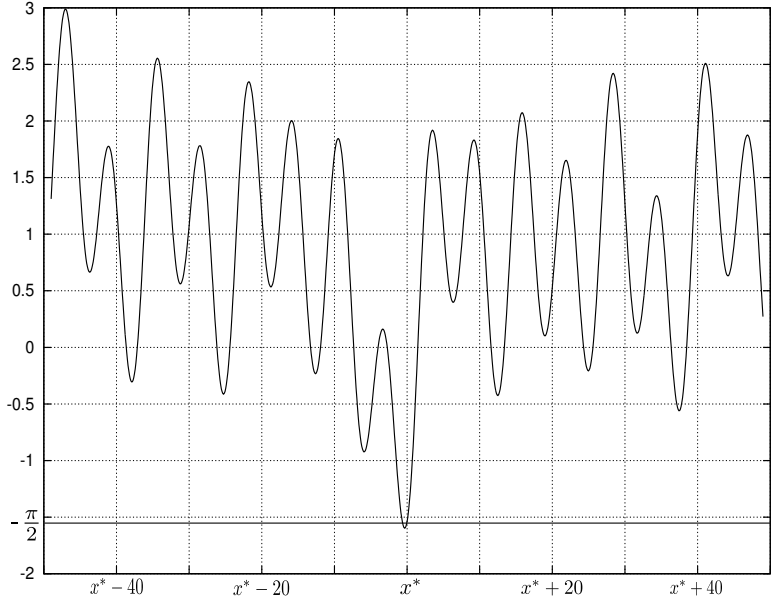} 
\caption{A region where $Q(x)<-\frac{\pi}{2}$. Here $x^*=8203872394818031742687.5 \times \pi$.} 
\label{fig2_xstar}
\end{figure}

As we see from Table \ref{tab2}, there are no values of $Q(\hat x^{(k)}_j)$ for $k\le 10$ which are smaller than $-\pi/2$. And $k=10$ was the largest value that we could possibly check on our desktop computer in a reasonable amount of time: in order to perform the computations for the next value $k=13$ (for which we need to evaluate 
$Q(\hat x^{(13)}_j)$, $1\le j \le K(13)=64411251125$) we would need to wait for nine months. 
 Luckily enough, after just two days of computations at the level $k=13$, we have found that at the index $j^*=63683535496$
we have
\beq\label{Q_xstar0}
Q(x^*)=-1.5398404\dots,
\eeq
where 
\beqq
x^*:=\hat x^{(13)}_{j^*}=(4j^*+3)K(13) \frac{\pi}{2}=8203872394818031742687.5 \times \pi.
\eeqq
The value in \eqref{Q_xstar0} is just slightly above $-\pi/2$. However, looking at the results in Table \ref{tab1}, we see that the local minimums of $Q(x)$ tend to occur near points $x$ such that $x/\pi = 0.4 \; ({\textnormal{mod}} \; 1)$. We check the value of $Q(x)$ at this point, and find that
\beq\label{Q_xstar}
Q(x^*-0.1\pi)=-1.5970415\dots.
\eeq

The value of $Q(x^*-0.1\pi)$ in \eqref{Q_xstar} was computed using the second approximation $Q_2(x)$, provided 
by Theorem \ref{thm_big_Oh_of_x_third}. 
Given the error estimate $O(x^{-\frac{1}{3}}\ln(x))$ in Theorem \ref{thm_big_Oh_of_x_third} and 
our previous discussion of the numerical results shown in Figure \ref{fig1_error}, 
we would expect that $Q_2(x)-Q(x) \approx 10^{-7}$ when $x \approx 10^{22}$. In order to confirm this and to check the accuracy of computation in \eqref{Q_xstar},  we have evaluated $Q_1(x^*-0.1\pi)$, which is the approximation given by Theorem \ref{thm_big_Oh_of_x_half} with $\epsilon=0.025$ and $M=4$. We found that the difference between the two results is less than $10^{-7}$, as expected. 
Finally, in order to eliminate the possibility that there was some loss of accuracy in \eqref{Q_xstar} due to the
addition of a large number of terms in \eqref{eqn_Qx_third}, we have evaluated $Q_2(x^*-0.1\pi)$ 
with the working precision of 60 decimal digits. We have used D. Bailey's MPFUN multi-precision library for Fortran90. Again, the results were in perfect agreement with \eqref{Q_xstar}. While we cannot prove
\eqref{Q_xstar} rigorously (since we do not know the explicit values for the implied constants in \eqref{eqn_Qx_third}), these numerical results
provide compelling evidence that 
\eqref{Q_xstar} is correct, and therefore, provides an explicit example for 
the key step of the disproof of the Clark and Ismail conjecture by Alzer et. al. \cite{Alzer2005}.

The plot of $Q(x)$ near the point $x^*$ is presented on Figure \eqref{fig2_xstar}. 

In conclusion, we would like to discuss the computation time of $Q(x^*-0.1\pi)$. This will help us to put the three approximations to $Q(x)$ into perspective. 
 It takes just five seconds to establish \eqref{Q_xstar} using the $O(x^{\frac{1}{3}})$ algorithm provided by Theorem \ref{thm_big_Oh_of_x_third}. Computing the same value 
using the $O(x^{\frac{1}{2}})$ algorithm presented in Theorem \ref{thm_big_Oh_of_x_half} took almost fourteen hours.
 By extrapolation, we find that the $O(x)$ algorithm from Proposition \ref{prop_big_Oh_of_x} would need more than forty million years in order to complete the same task.

\section{Proofs}\label{section_proofs}

\vspace{0.25cm}
{\it Proof of Theorem \ref{thm_big_Oh_of_x_half}:}
We set $h:=1/x$ and apply the Euler-Maclaurin summation formula
to function $G(x)=\frac{1}{x} \sin\left(\frac{1}{x}\right)$: for every $j\ge 0$
\beq\label{EM_formula}
h \left[\frac{1}{2} G(Nh)+\sum\limits_{n\ge N+1} 
G(nh) \right]&=&\int\limits_{Nh}^{\infty} G(u) \d u - \sum\limits_{k=1}^{j}  \frac{B_{2k}}{(2k)!} h^{2k} G^{(2k-1)}(Nh) \\ \nonumber
&+&\frac{h^{2j+1}}{(2j+1)!}\int\limits_{Nh}^{\infty} {\textnormal{B}}_{2j+1}(\{u/h\}) G^{(2j+1)}(u) \d u,
\eeq
where ${\textnormal{B}}_n(\cdot)$ denote Bernoulli polynomials and $B_n={\textnormal{B}}_n(0)$ are Bernoulli numbers.

Changing the variable of integration $u=1/v$ and using \eqref{def_Si} we find that the first integral in the right-hand side of \eqref{EM_formula} is equal to 
$\Si(x/N)$. 

Let us now estimate the second integral in the right-hand side of \eqref{EM_formula}.  Using induction it is easy to prove that there exists a sequence of polynomials $\{p_{n}(u)\}_{n\ge 0}$, such that ${\textnormal{deg}}(p_n)=n$ and for each $j\ge 0$ we have
\beq\label{G_derivative}
G^{(2j+1)}(u)=u^{-2j-2} \left[ \cos(u^{-1}) p_{2j+1}(u^{-1} )+\sin(u^{-1}) p_{2j}(u^{-1}) \right].
\eeq
Formula \eqref{G_derivative} implies that for every $j \ge 0$ there exists a constant $C_1(j)>0$ such that for all $u>0$ we have
\beq\label{bound_G_derivatives}
|G^{(2j+1)}(u)| < C_1(j)\left( u^{-4j-3}+u^{-2} \right)
\eeq
The periodic function  $|{\textnormal{B}}_{2j+1}(\{u/h\})|$ is bounded from above by some positive constant $C_2(j)$, therefore we obtain the estimate
\beqq
\left | \; \int\limits_{Nh}^{\infty} {\textnormal{B}}_{2j+1}(\{u/h\}) G^{(2j+1)}(u) \d u \right| 
< C_1(j) \times C_2(j) \times  \int\limits_{Nh}^{\infty} \left( u^{-4j-3}+u^{-2} \right) \d u= O_j\left((Nh)^{-4j-2}\right).
\eeqq
Using the fact that $Nh\approx x^{-\frac{1}{2}+\epsilon}$ as $x\to +\infty$, we conclude
\beq\label{int2_estimate}
\frac{h^{2j+1}}{(2j+1)!}\int\limits_{Nh}^{\infty} {\textnormal{B}}_{2j+1}(\{u/h\}) G^{(2j+1)}(u) \d u=O_j\left(x^{-(4j+2)\epsilon} \right).
\eeq

Combining  \eqref{Q_main_sum_1}, \eqref{EM_formula} and \eqref{int2_estimate} we see that for every $j\ge 0$
\beq\label{E_N_x_estimate1}
{\mathcal E}_N(x) =-\sum\limits_{k=1}^{j}  \frac{B_{2k}}{(2k)!} 
\frac{G^{(2k-1)}\left( \frac{N}{x} \right)}{x^{2k}} 
 +O_j\left(x^{-(4j+2)\epsilon} \right).
\eeq
For every $\epsilon \in (0,\frac12)$ and $M \ge 0$ we can find $j$ large enough so that 
$(4j+2)\epsilon> \frac{1}{2} + (4M+3) \epsilon$. Also, from \eqref{bound_G_derivatives} we obtain
\beqq
\frac{G^{(2k-1)}\left( \frac{N}{x} \right)}{x^{2k}} = x^{-2k} O_k( (N/x)^{-4k+1} )= O_{k}(x^{-\frac{1}{2}-(4k-1)\epsilon}).
\eeqq
Combining these two facts with \eqref{E_N_x_estimate1} we conclude that \eqref{error_term} holds. 
\qed

The proof of the Theorem \ref{thm_big_Oh_of_x_third} is based on the following two Lemmas. 
 
\begin{lemma}\label{Lemma_12_derivative_test}
\emph{(First and Second derivative test, see Lemmas 5.1.2 and 5.1.3 in \cite{Huxley1996})}

\noindent
Let $f: (\alpha,\beta) \mapsto \r$ be twice continuously differentiable and 
$g: (\alpha, \beta) \mapsto [0,\infty)$ be monotone. For $i \in \{1,2\}$ we define
$\lambda_i=\inf \{ |f^{(i)}(u)| \; : \; \alpha < u < \beta \}$ and assume that $\lambda_2>0$. 
Then 
\beq
\left | \int\limits_{\alpha}^{\beta} g(u) \sin(f(u)) \d u \right| < 16 \times \max \{ g(\alpha+), g(\beta-)\}
\times \min\left\{ \frac{1}{\lambda_1}, \frac{1}{\sqrt{\lambda_2} } \right\}
\eeq
\end{lemma}

\begin{lemma}\label{Lemma_Poisson_sum_formula}
\emph{(Truncated Poisson summation formula, see Lemma 5.4.3 in \cite{Huxley1996} or Lemma 4.10 in 
\cite{Titchmarsh1986})}

\noindent
Let $f: [a,b] \mapsto \r$ have a continuous and decreasing derivative $f'(u)$. Let $g: [a,b] \mapsto [0,\infty)$ 
be decreasing function, such that $|g'(u)|$ is also decreasing. Define $\alpha=f'(b)$ and 
$\beta =f'(a)$. Then
\beq\label{eqn_Poisson_sum}
\sum\limits_{a< n \le b} g(n) e^{2\pi i f(n)} = 
\sum\limits_{\alpha-\frac{1}{4} < m < \beta + \frac{1}{4} } 
\int\limits_{a}^{b} g(u) e^{2\pi \i ( f(u)- m u)} \d u + O(g(a) \ln(\beta-\alpha+2))+O(|g'(a)|)
\eeq
\end{lemma}

\vspace{0.25cm}
\noindent
{\it Proof of Theorem \ref{thm_big_Oh_of_x_third}:}
The proof will be carried in three steps. Our first goal is to  establish the following identity: 
\beq\label{proof_thm2_identity1}
\sum\limits_{n\ge N} \frac{1}{n} \sin\left(\frac{x}{n}\right)=
\sum\limits_{0\le m < M+\frac{1}{4}} \left[ \;\int\limits_N^{\infty} \frac{1}{u} \sin\left(\frac{x}{u}+2 \pi  m  u \right) \d u \right] + 
O( N^{-1} \ln(x)). 
\eeq
Note that the conditions $2\pi M N^2=x$, $M\ge 1$, $N\ge 1$ imply that
\beq\label{main_conditions}
M=O(x), \;\;\; N=O(\sqrt{x}), \;\;\;  x \ge 2\pi. 
\eeq
We define $a:=N$, $b:=x^2$,  $f(u):=-\frac{x}{2\pi u}$ and $g(u):=\frac{1}{u}$ and check that all conditions of Lemma \ref{Lemma_Poisson_sum_formula} are satisfied. We find that $\alpha=f'(b)=\frac{1}{2\pi x^3}$ and $\beta=f'(a)=M$
 and note that for $x \ge 2 \pi$ we have $-1<\alpha -\frac{1}{4} <0$. 
 Applying  Lemma \ref{Lemma_Poisson_sum_formula} and taking imaginary part of both sides of \eqref{eqn_Poisson_sum}
 we obtain
\beq\label{proof_thm2_identity11}
\sum\limits_{N < n \le x^2} \frac{1}{n} \sin\left(\frac{x}{n} \right)=
\sum\limits_{0\le m < M+\frac{1}{4}} \left[ \; \int\limits_N^{x^2} \frac{1}{u} 
 \sin\left(\frac{x}{u}+2 \pi  m  u \right)  \d u \right] +  O(N^{-1} \ln(M)). 
\eeq

Next, assume that $m\ge 1$ and $w>x^2$. We define $f(u):=\frac{x}{u}+2 \pi  m  u$ and note that 
\beqq
\lambda_1&=&\inf \{ |f'(u)| \; : \; x^2 < u < w \}=-\frac{1}{x^3}+2\pi m, \\
\lambda_2&=&\inf \{ |f''(u)| \; : \; x^2 < u < w \}=\frac{2}{x^5}>0.
\eeqq
We apply Lemma \ref{Lemma_12_derivative_test} (the first derivative test) with the above function $f(u)$ and $g(u):=u^{-1}$ and conclude that for $m\ge 1$, $x \ge 2\pi$ and any $w>x^2$ 
\beqq
\int\limits_{x^2}^{w} \frac{1}{u} \sin\left(\frac{x}{u}+2 \pi  m  u \right) \d u 
 \ll \frac{16}{x^2} \frac{1}{\left(-\frac{1}{x^3}+2\pi m\right)} \ll \frac{1}{x^2}.
\eeqq
Therefore for all $x \ge  2\pi$ and $m\ge 1$ we have
\beqq
\int\limits_{x^2}^{\infty} \frac{1}{u} \sin\left(\frac{x}{u}+2 \pi  m  u \right) \d u \ll \frac{1}{x^2}
\eeqq
When $m=0$ we have the following trivial estimate
\beqq
\int\limits_{x^2}^{\infty} \frac{1}{u} \sin\left(\frac{x}{u} \right) \d u 
\ll \int\limits_{x^2}^{\infty} \frac{1}{u} \left | \sin\left(\frac{x}{u} \right) \right| \d u  \ll \int\limits_{x^2}^{\infty} \frac{x}{u^2} \d u =\frac{1}{x}.
\eeqq
The above two estimates and the fact that $M=O(x)$  give us 
\beq\label{estimate_sum_integrals}
\sum\limits_{0\le m < M+\frac{1}{4}} \left[ \; \int\limits_{x^2}^{\infty} \frac{1}{u} 
 \sin\left(\frac{x}{u}+2 \pi  m  u \right)  \d u \right] \ll \frac{M}{x^{2}} + \frac{1}{x} \ll \frac{1}{x}. 
\eeq
Combining \eqref{main_conditions}, \eqref{proof_thm2_identity11}, \eqref{estimate_sum_integrals} and the following simple estimate
\beqq
 \sum\limits_{n \ge x^2 } \frac{1}{n} \sin\left(\frac{x}{n}\right) \ll 
 \sum\limits_{n \ge x^2 } \frac{x}{n^2} \ll \frac{1}{x}, 
\eeqq
we obtain \eqref{proof_thm2_identity1}.

The second step consists in establishing the following result
\beq\label{proof_thm2_identity2}
\sum\limits_{0\le m < M+\frac{1}{4}} \left[ \;\int\limits_0^{N} \frac{1}{u} \sin\left(\frac{x}{u}+2 \pi  m  u \right) \d u \right] \ll N^{-1} \ln(x)+x^{-\frac12}N^{\frac{1}{2}}.
\eeq
We change the variable of integration $u=1/v$ and obtain
\beq\label{integral_identity}
 \int\limits_{0}^{N} \frac{1}{u} \sin\left(\frac{x}{u}+2 \pi  m  u \right) \d u =
 \int\limits_{N^{-1}}^{\infty} \frac{1}{v} \sin\left(xv+ \frac{2 \pi  m}{v} \right) \d v.  
\eeq
Define  $f(v):=xv+ 2 \pi  m/v$ and assume that $0\le m  \le \hat M$, where $\hat M:=\lfloor M-\frac{1}{4} \rfloor$. 
Then for $v > N^{-1}$ we have
\beqq 
f'(v)=x-\frac{2\pi m}{v^2}>x-2\pi m N^2. 
\eeqq 
Note that $ m  \le \hat M \le M-\frac{1}{4}$ implies
\beqq
x-2\pi m N^2\ge x-2\pi \left(M-\frac{1}{4}\right) N^2=\frac{\pi}{2N^2},
\eeqq
thus $\lambda_1=\inf \{ |f'(v)| \; : \; N^{-1} < v < \infty \}=x-2\pi m N^2$ and $\lambda_1>0$. 
Applying Lemma \ref{Lemma_12_derivative_test} (the first derivative test) to the integral in the right-hand side of \eqref{integral_identity}
we conclude that
\beqq
 \int\limits_{0}^{N} \frac{1}{u} \sin\left(\frac{x}{u}+2 \pi  m  u \right) \d u =
 \int\limits_{N^{-1}}^{\infty} \frac{1}{v} \sin\left(xv+ \frac{2 \pi  m}{v} \right) \d v  \ll 
 \frac{16N}{x-2\pi m N^2}. 
\eeqq
Using the above estimate we obtain
\beq\label{sum_m_0_hat_M}
&&\sum\limits_{m=0}^{\hat M}
\left[ \; \int\limits_{0}^{N} \frac{1}{u} \sin\left(\frac{x}{u}+2 \pi  m  u \right) \d u  \right] \ll 
\sum\limits_{m=0}^{\hat M} \frac{16N}{x-2\pi m N^2}=
\sum\limits_{j=0}^{\hat M} \frac{16N}{x-2\pi (\hat M - j) N^2} \\ \nonumber
&\le&  \sum\limits_{j=0}^{\hat M} 
\frac{16 N}{x-2\pi (M-\frac{1}{4}) N^2 +2\pi  j N^2}=
\sum\limits_{j=0}^{\hat M} 
\frac{32 }{\pi N(1+4j)} \ll N^{-1} \ln(2+\hat M) \ll N^{-1} \ln(x),
\eeq
where we have changed the index of summation $j=\hat M -m$ and have used \eqref{main_conditions} and the fact that 
$\sum_{1\le n \le y} (1/n) \ll \ln(y)$. 

In order to finish the proof of \eqref{proof_thm2_identity2} we have to  consider the possible
integral term with $\hat M < m < M+\frac{1}{4}$ in the sum in the left-hand side of 
\eqref{proof_thm2_identity2}. We again change the variable of integration $v=1/u$ and consider the following two integrals
\beq\label{term_close_to_M}
\int\limits_{0}^{N} \frac{1}{u} \sin\left(\frac{x}{u}+2 \pi  m  u \right) \d u 
 =\int\limits_{N^{-1}}^{2N^{-1}} \frac{1}{v} \sin\left(xv+ \frac{2 \pi  m}{v} \right) \d v 
 +\int\limits_{2N^{-1}}^{\infty} \frac{1}{v} \sin\left(xv+ \frac{2 \pi  m}{v} \right) \d v 
\eeq
Define $f(u):=xu + \frac{2\pi m}{u}$, then $f'(u)=x- \frac{2\pi m}{u^2}$ and 
$f''(u)=\frac{4\pi m}{u^3}$. It is easy to see that 
\beqq
\lambda_2=\min \{|f''(u)| \; : \; N^{-1} \le u \le 2 N^{-1} \}=\frac{1}{2}\pi m N^3.
\eeqq
Applying Lemma \ref{Lemma_12_derivative_test} (the second derivative test) we obtain the following estimate
\beq\label{term_close_to_M1}
\int\limits_{N^{-1}}^{2N^{-1}} \frac{1}{u} \sin\left(xu+ \frac{2 \pi  m}{u} \right) \d u \ll (mN)^{-\frac{1}{2}} \ll x^{-\frac12}N^{\frac{1}{2}}
\eeq
where in the last step we have used the facts that $m> M-\frac{3}{4}$ and $MN=x/(2\pi N)$. 
Next, we find that for $u>2N^{-1}$ 
\beqq
f'(u)=x- \frac{2\pi m}{u^2}>x-\frac{1}{2}\pi (M+\frac14) N^2 = 
\frac{3x}{4}-\frac{\pi}{8} N^2 >\frac{x}{2},
\eeqq
thus by Lemma \ref{Lemma_12_derivative_test} (the first derivative test) we obtain
\beq\label{term_close_to_M2}
\int\limits_{2N^{-1}}^{\infty} \frac{1}{u} \sin\left(xu+ \frac{2 \pi  m}{u} \right) \d u \ll Nx^{-1} \ll N^{-1},
\eeq
where in the last step we have used the upper bound $N x^{-1}=N^{-1}/(2\pi M) \ll N^{-1}$. 
Equations \eqref{term_close_to_M}, \eqref{term_close_to_M1} and \eqref{term_close_to_M2} show that 
for $\hat M < m < M+\frac{1}{4}$
\beqq
\int\limits_{0}^{N} \frac{1}{u} \sin\left(\frac{x}{u}+2 \pi  m  u \right) \d u \ll N^{-1} + x^{-\frac{1}{2}} N^{\frac{1}{2}}. 
\eeqq
 The above estimate  combined with 
\eqref{sum_m_0_hat_M} complete the proof of \eqref{proof_thm2_identity2}. 

Now we have all the ingredients for the last step of the proof of Theorem \ref{thm_big_Oh_of_x_third}. 
We combine the two results \eqref{proof_thm2_identity1} and \eqref{proof_thm2_identity2} and conclude that
\beq\label{last_identity}
\sum\limits_{n\ge N} \frac{1}{n} \sin\left(\frac{x}{n}\right)&=&
\sum\limits_{0\le m < M+\frac{1}{4}} \left[ \;\int\limits_0^{\infty} \frac{1}{u} \sin\left(\frac{x}{u}+2 \pi  m  u \right) \d u \right] + 
O( N^{-1} \ln(x)+x^{-\frac12}N^{\frac{1}{2}})\\ \nonumber
&=&\frac{\pi}{2}+
\pi \sum\limits_{1\le m \le M-\frac{1}{4}}  J_0(2\sqrt{2 \pi m x} ) + O( N^{-1} \ln(x)+x^{-\frac12}N^{\frac{1}{2}})
\eeq
where we have used the following integral identities (see formulas 3.721.1 and 3.868.1 in \cite{Jeffrey2007}):
\beqq
\int\limits_{0}^{\infty} \frac{1}{u} \sin\left(\frac{x}{u} \right) \d u= \frac{\pi}2, \;\;\;
\int\limits_0^{\infty} \frac{1}{u} \sin\left(\frac{x}{u}+2 \pi  m  u \right) \d u=
\pi J_0(2\sqrt{2 \pi m x} ). 
\eeqq
In order to obtain \eqref{eqn_Qx_third} one has to combine \eqref{last_identity}, the asymptotic expression 
\eqref{J_0_asymptotics} for $J_0(2\sqrt{2 \pi m x} )$ and 
the following simple estimate 
\beqq
\sum\limits_{1 \le m < M +\frac{1}{4}} (xm)^{-\frac{3}{4}} \ll 
x^{-\frac{3}{4}} M^{\frac{1}{4}} \ll x^{-\frac{1}{2}} \ll N^{-1}.
\eeqq
\qed

The proof of Theorem \ref{thm_big_Oh_of_x_eps} is based on the following result from the theory of exponential sums.
\begin{theorem}\label{thm29}
\emph{(Theorem 2.9 in \cite{Graham1991})}
Let $q\ge 0$ be an integer and $Q=2^q$. Suppose that $f$ has $q+2$ continuous derivatives on an interval $I$ and that $I \subset (N,2N]$. Assume also that there is some constant $F$ such that 
for $1\le r \le q+2$ we have
\beq\label{estimate_derivatives}
|f^{(r)}(u)| \approx F N^{-r}, \;\;\; u\in I.
\eeq
 Then
\beq\label{thm29_main_estimate}
\sum\limits_{n \in I} e^{2\pi \i f(n)} \ll F^{\frac{1}{4Q-2}} N^{1-\frac{q+2}{4Q-2}}+F^{-1} N.
\eeq
The implied constant in \eqref{thm29_main_estimate} depends only upon the implied constants in  
\eqref{estimate_derivatives}. 
\end{theorem}

\vspace{0.25cm}
\noindent
{\it Proof of Theorem \ref{thm_big_Oh_of_x_eps}:}
Let $N$ be any  positive number. For $u \ge N$ we define
\beqq
S(u):=\sum\limits_{N \le n \le u} \sin\left(\frac{x}{n} \right)
\eeqq
and $S(u):=0$ for $u<N$. 
We also define $f(u):=x/(2\pi u)$ and $F:=x/N$. Computing $f^{(r)}(u)$ one can easily check  that for $N \le u \le 2N$ and $1\le r \le q+2$ 
\beqq
\frac{r!}{2^{2+r} \pi}  \times FN^{-r} \le |f^{(r)}(u)| \le \frac{r!}{2 \pi} \times  FN^{-r},
\eeqq
therefore all
the conditions of Theorem \ref{thm29} are satisfied. Taking the imaginary part of \eqref{thm29_main_estimate} we conclude that for all $q\ge 0$, $N\ge 1$ and $N \le u \le 2N$ 
\beq\label{estimate_Su}
S(u)\ll_q  x^{\frac{1}{4Q-2}}N^{1-\frac{q+3}{4Q-2}} + x^{-1} N^2 
\eeq
where $Q=2^q$. We apply the integration by parts and find that for all $N\le M < 2N$ we have
\beq\label{thm3_proof_estimate1}
\sum\limits_{n=N}^{M} \frac{1}{n} \sin\left(\frac{x}{n} \right)
&=&\int\limits_{N-\frac{1}{2}}^{M+\frac{1}{2}} \frac{\d S(u)}{u} =
\frac{S(M+\frac{1}{2})}{M+\frac{1}{2}} +\int\limits_{N-\frac{1}{2}}^{M+\frac{1}{2}} 
\frac{S(u)}{u^2}  \d u \\ \nonumber 
&\ll& N^{-1} \max\left\{|S(u)| \; : N-\frac{1}{2} \le u \le M+\frac{1}{2} \right\}
\ll_q x^{\frac{1}{4Q-2}}N^{-\frac{q+3}{4Q-2}} + x^{-1} N,
\eeq
where in the last step we have used \eqref{estimate_Su}. 

Let us define $N:=\lfloor x^{\epsilon} \rfloor$  and
\beqq
 j=j(x):=\lfloor\log_2(N^{-1}x^{3/5} )\rfloor, \;\;\;\;\;\;
\alpha=\alpha(q,\epsilon):=((q+3)\epsilon-1)/(4Q-2).
\eeqq 
Dividing the interval $N \le n < x^{3/5}$ into subintervals of the form $N_1 \le n \le N_2 < 2N_1$ we obtain
\beq\label{thm3_proof_estimate2}
\sum\limits_{N \le n < x^{3/5}}  \frac{1}{n} \sin\left(\frac{x}{n} \right)
&=&\sum\limits_{i=0}^{j-1} \sum\limits_{n=2^i N}^{2^{i+1} N-1} \frac{1}{n} \sin\left(\frac{x}{n} \right)
+\sum\limits_{2^j N \le n < x^{3/5}} \frac{1}{n} \sin\left(\frac{x}{n} \right) \\ \nonumber
&\ll_q& 
x^{\frac{1}{4Q-2}}N^{-\frac{q+3}{4Q-2}} (j+1) + x^{-1}(x^{3/5}-N) 
\ll_q x^{-\alpha}\ln(x) + x^{-\frac{2}{5}},
\eeq
where we have used the facts that $j\ll \ln(x)$ and $N \approx x^{\epsilon}$. 
 We choose $q=\lfloor \epsilon^{-1} \rfloor -1$ and find that
\beqq
\alpha=\frac{ \lfloor \epsilon^{-1} \rfloor \epsilon + 2 \epsilon-1}{2^{\lfloor \epsilon^{-1} \rfloor +1}-2} >
 \frac{\epsilon}{2^{\lfloor \epsilon^{-1} \rfloor +1}-2}>{\epsilon}2^{-\frac{1}{\epsilon}}=\delta(\epsilon).
\eeqq
This fact combined with the estimate \eqref{thm3_proof_estimate2} shows that
\beq\label{thm3_final_estimate}
\sum\limits_{x^{\epsilon} \le n < x^{3/5}}  \frac{1}{n} \sin\left(\frac{x}{n} \right)=O_{\epsilon}(x^{-\delta}).
\eeq
Taking $\epsilon=1/10$ in \eqref{Qx_truncated_sum} we find that
\beqq
Q(x)=\sum\limits_{1 \le n < x^{3/5}} \frac{1}{n} \sin\left(\frac{x}{n}\right)
+\frac{\pi}2 + O(x^{-\frac{2}{5}}).
\eeqq
Combining this result with \eqref{thm3_final_estimate} ends the proof of Theorem \ref{thm_big_Oh_of_x_eps}.
\qed



\begin{thebibliography}{10}

\bibitem{Alzer2005}
H.~Alzer, C.~Berg, and S.~Koumandos.
\newblock On a conjecture of {C}lark and {I}smail.
\newblock {\em J. Approx. Theory}, 134(1):102--113, 2005.

\bibitem{Clark2003}
{\relax W.E}.~Clark and {\relax M.E.H}.~Ismail.
\newblock Inequalities involving gamma and psi functions.
\newblock {\em Anal. Appl.}, 1(1):129--140, 2003.

\bibitem{Codeca1984}
P.~Codec\`a.
\newblock On the properties of oscillation and almost periodicity of certain
  convolutions.
\newblock {\em Rend. Sem. Mat. Univ. Padova}, 71:103--119, 1984.

\bibitem{Delange1980}
H.~Delange.
\newblock Sur la fonction $f(x)=\sum\limits_{n=1}^{\infty} (1/n) \sin(x/n)$.
\newblock {\em Theorie analytique et elementaire des nombres, Caen, 29-30
  Septembre 1980, Journees mathematiques SMF-CNRS}, 1980.

\bibitem{Doring1985}
B.~D{\"o}ring.
\newblock On the zeros of {F}lett's function.
\newblock {\em J. Comput. Appl. Math.}, 12 - 13:265 -- 270, 1985.

\bibitem{Flett1950}
{\relax T.M}.~Flett.
\newblock On the function $\sum\limits_{n=1}^{\infty} (1/n) \sin(t/n)$.
\newblock {\em J. London Math. Soc.}, 25:5--19, 1950.

\bibitem{Gautschi2005}
W.~Gautschi.
\newblock The {H}ardy-{L}ittlewood function: an exercise in slowly convergent
  series.
\newblock {\em J. Comput. Appl. Math.}, 179(1-2):249--254, 2005.

\bibitem{Jeffrey2007}
I.~S. Gradshteyn and I.~M. Ryzhik.
\newblock {\em Table of integrals, series, and products}.
\newblock Elsevier/Academic Press, Amsterdam, seventh edition, 2007.

\bibitem{Graham1991}
S.~Graham and G.~Kolesnik.
\newblock {\em Van der Corput's method of exponential sums}, volume 126 of {\em
  London Mathematical Society Lecture Notes Series}.
\newblock Cambridge University Press, 1991.

\bibitem{HL1936}
{\relax G.H}.~Hardy and {\relax J.E}.~Littlewood.
\newblock Notes on the theory of series ({XX}): on {L}ambert series.
\newblock {\em Proc. London Math. Soc.}, s2-41(1):257--270, 1936.

\bibitem{Huxley1996}
M.~Huxley.
\newblock {\em Area, Lattice Points and Exponential Sums}.
\newblock Oxford University Press, 1996.

\bibitem{Petermann1988}
Y.-F. P\'etermann.
\newblock About a theorem of {P}aolo {C}odec\`a's and {$\Omega$}-estimates for
  arithmetical convolutions.
\newblock {\em Journal of Number Theory}, 30(1):71 -- 85, 1988.

\bibitem{Segal1972}
{\relax S.L}.~Segal.
\newblock On $\sum (1/n) \sin(x/n)$.
\newblock {\em J. London Math. Soc.}, s2-4(3):385--393, 1972.

\bibitem{Titchmarsh1986}
E.~Titchmarsh.
\newblock {\em The theory of the Riemann zeta-function}.
\newblock Oxford University Press, second edition, 1986.

\bibitem{Lune1981}
J.~{\relax van de Lune}.
\newblock A note on the zeros of {F}lett's function.
\newblock {\em Afdeling Zuivere Wiskunde, Report ZW 167, Mathematisch Centrum,
  Amsterdam}, 1981.

\bibitem{Walfisz1936}
A.~Walfisz.
\newblock {\em Weilsche {E}xponentialsummen in der neueren {Z}ahlentheorie}.
\newblock Math Forschugsber, 15, V.E.B. Deutcher Verlag der Wiss., 1936.

\end{thebibliography}

\end{document}